\documentclass[12pt]{amsart}
\usepackage[colorlinks=true, pdfstartview=FitV, linkcolor=blue, citecolor=blue, urlcolor=blue]{hyperref}

\usepackage{amssymb,amsmath, amscd}
\usepackage{times, verbatim}
\usepackage{graphicx}
 \usepackage[usenames, dvipsnames]{color}
\usepackage{amsmath,amssymb,amsfonts}
\usepackage{enumerate}

\usepackage{tikz}
\usetikzlibrary{calc, arrows, decorations.markings} \tikzset{>=latex}
\usepackage{amsmath,amssymb,amsfonts}
\usepackage{amscd, amssymb, latexsym, amsmath, amscd}
\usepackage[all]{xy}
\usepackage{pb-diagram}
\usepackage{enumerate}
\usepackage{ulem}
\usepackage{anysize}
\marginsize{3cm}{3cm}{3cm}{3cm}
\usepackage[colorlinks=true, pdfstartview=FitV, linkcolor=blue, citecolor=blue, urlcolor=blue]{hyperref}
\usepackage{multido,pstricks}
\usepackage{amssymb,amsmath, amscd}
\usepackage{times, verbatim}
\usepackage{graphicx}
\usepackage[english]{babel}
\usepackage{cleveref}
\usepackage{mathrsfs}
\usepackage{tikz}
\usetikzlibrary{arrows.meta}
\usetikzlibrary{bending}
\usetikzlibrary{calc, arrows, decorations.markings} \tikzset{>=latex}
\usepackage{amsmath,amssymb,amsfonts}
\usepackage{amscd, amssymb, latexsym, amsmath, amscd}
\usepackage[all]{xy}
\usepackage{parskip}
\usepackage{pb-diagram}
\usepackage{enumerate}
\usepackage{ulem}
\usepackage{anysize}
\marginsize{3cm}{3cm}{3cm}{3cm}
\input xy
\xyoption{all}
\usepackage{pb-diagram}
\usepackage[all]{xy}
\usepackage{tabularx}
\usepackage{ytableau}
\usepackage{xcolor}
\input xy
\xyoption{all}

\begin{document}
\theoremstyle{plain}

\newtheorem{theorem}{Theorem}[section]
\newtheorem{thm}[equation]{Theorem}
\newtheorem{prop}[equation]{Proposition}
\newtheorem{cor}[equation]{Corollary}
\newtheorem{conj}{Conjecture}
\newtheorem{lemma}[equation]{Lemma}
\newtheorem{define}{Definition}
\newtheorem{question}[equation]{Question}

\theoremstyle{definition}
\newtheorem{conjecture}[theorem]{Conjecture}

\newtheorem{example}{Example}
\numberwithin{equation}{section}

\newtheorem{remark}{Remark}

\newcommand{\Hecke}{\mathcal{H}}
\newcommand{\Liea}{\mathfrak{a}}
\newcommand{\Cmg}{C_{\mathrm{mg}}}
\newcommand{\Cinftyumg}{C^{\infty}_{\mathrm{umg}}}
\newcommand{\Cfd}{C_{\mathrm{fd}}}
\newcommand{\Cinftyfd}{C^{\infty}_{\mathrm{ufd}}}

\newcommand{\sspace}{\Gamma \backslash G}
\newcommand{\PP}{\mathcal{P}}
\newcommand{\bfP}{\mathbf{P}}
\newcommand{\bfQ}{\mathbf{Q}}
\newcommand{\Siegel}{\mathfrak{S}}
\newcommand{\g}{\mathfrak{g}}
\newcommand{\A}{\mathcal{A}}
\newcommand{\B}{\mathcal{B}}
\newcommand{\Q}{\mathbb{Q}}
\newcommand{\Gm}{\mathbb{G}_m}
\newcommand{\Nm}{\mathbb{N}m}
\newcommand{\ii}{\mathfrak{i}}
\newcommand{\II}{\mathfrak{I}}

\newcommand{\kk}{\mathfrak{k}}
\newcommand{\nn}{\mathfrak{n}}
\newcommand{\tF}{\tilde{F}}
\newcommand{\p}{\mathfrak{p}}
\newcommand{\m}{\mathfrak{m}}
\newcommand{\bb}{\mathfrak{b}}
\newcommand{\Ad}{{\rm Ad}\,}
\newcommand{\ttt}{\mathfrak{t}}
\newcommand{\frakt}{\mathfrak{t}}
\newcommand{\U}{\mathcal{U}}
\newcommand{\Z}{\mathbb{Z}}

\newcommand{\LL}{{\mathcal L}}
\newcommand{\bfG}{\mathbf{G}}
\newcommand{\bfT}{\mathbf{T}}
\newcommand{\R}{\mathbb{R}}
\newcommand{\ST}{\mathbb{S}}
\newcommand{\h}{\mathfrak{h}}
\newcommand{\bC}{\mathbb{C}}
\newcommand{\C}{\mathbb{C}}
\newcommand{\F}{\mathbb{F}}
\newcommand{\N}{\mathbb{N}}
\newcommand{\qH}{\mathbb {H}}
\newcommand{\temp}{{\rm temp}}
\newcommand{\Hom}{{\rm Hom}}
\newcommand{\Aut}{{\rm Aut}}
\newcommand{\rk}{{\rm rk}}
\newcommand{\Ext}{{\rm Ext}}
\newcommand{\End}{{\rm End}\,}
\newcommand{\Ind}{{\rm Ind}}
\newcommand{\ind}{{\rm ind}}
\newcommand{\Irr}{{\rm Irr}}
\def\circG{{\,^\circ G}}
\def\M{{\rm M}}
\def\diag{{\rm diag}}
\def\Ad{{\rm Ad}}
\def\As{{\rm As}}
\def\wG{{\widehat G}}
\def\wT{{\widehat T}}
\def\wB{{\widehat B}}
\def\wM{{\widehat M}}
\def\wlambda{{\widehat \lambda}}
\def\walpha{{\widehat \alpha}}
\def\wbeta{{\widehat \beta}}
\def\G{{\rm G}}
\def\H{{\rm H}}
\def\T{{\rm T}}
\def\SL{{\rm SL}}
\def\PSL{{\rm PSL}}
\def\GSp{{\rm GSp}}
\def\PGSp{{\rm PGSp}}
\def\Sp{{\rm Sp}}
\def\St{{\rm St}}
\def\GU{{\rm GU}}
\def\SU{{\rm SU}}
\def\U{{\rm U}}
\def\GO{{\rm GO}}
\def\GL{{\rm GL}}
\def\PGL{{\rm PGL}}
\def\GSO{{\rm GSO}}
\def\GSpin{{\rm GSpin}}
\def\GSp{{\rm GSp}}

\def\Gal{{\rm Gal}}
\def\SO{{\rm SO}}
\def\O{{\rm  O}}
\def\Sym{{\rm Sym}}
\def\sym{{\rm sym}}
\def\St{{\rm St}}
\def\Sp{{\rm Sp}}
\def\tr{{\rm tr\,}}
\def\ad{{\rm ad\, }}
\def\Ad{{\rm Ad\, }}
\def\rank{{\rm rank\,}}


\title{Some questions in Diophantine approximation, \\ real and p-adics}

\begin{abstract}
   The {\it Weak approximation theorem} describes the closure of $\G(\Q)$ inside $\G(\Q_p)$ as well as inside
  $\G(\R)$ for $\G$ an algebraic group over $\Q$; the closure is  always an open normal subgroup with finite abelian quotient, and is well understood in a certain sense even if precise
  results are not always available (such as for tori!). In this paper,
  for a finitely generated subgroup ${\mathcal L}\subset \G(\Q)$ we consider the topological closure of ${\mathcal L}$ inside $\G(\Q_p)$
  and $\G(\R)$. The paper is written mostly for $\G$ a torus or an abelian variety, but eventually considers  a variant of the question for 
  $\G$ a semisimple group. The paper is written with the wishful thinking
  that when dealing with questions on topological closure of {\it algebraic points} in an algebraic group defined over a number field,
  the simplest answers hold, a well-known
  principle known as ``Occum's razor''.
\end{abstract}

\author{Dipendra Prasad}
\address{Indian Institute of Technology Bombay, Powai, Mumbai-400076}
\email{prasad.dipendra@gmail.com}
\date{\today}

\thanks{These are the notes of some of the lectures I have given in the last few years,
  for which I have decided to keep the exposition, and also a bit of
  informality.
  I thank  ANRF, India for its support
through the JC Bose
National Fellowship of the Govt. of India, project number JBR/2020/000006.}

\maketitle
    
    \tableofcontents

    \newpage
    
    \section{Introduction}
    One of the first and most famous theorem on Diophantine approximations
    is the following theorem of Kronecker:
    
    \begin{thm}(Kronecker) A real number $\theta$, considered as an element
      of $\R/\Z$,  generates a dense subgroup of $\R/\Z$ if and only if $\theta$ does not
      belong to $\Q$.
    \end{thm}

    This theorem has many refinements already in one variable of great significance, such as the following theorem of Dirichlet as a consequence of the
    ``pigeonhole principle'':

    \begin{thm}(Dirichlet)
      Given $\theta \in \R$, an irrational number,
      \[ \left | \theta - \frac{p}{q} \right | < \frac{1}{q^2},\]
      for infinitely many rational numbers $\frac{p}{q}$.
    \end{thm}

    On the other hand,
    a theorem of      Thue–Siegel–Roth says that we cannot do much better:
 \begin{thm} (Thue–Siegel–Roth)
   Given $\theta$ algebraic and irrational, for any $\epsilon >0$,
   there exists a constant $C(\theta, \epsilon) >0$ such that:
      \[ \left | \theta - \frac{p}{q} \right | >  \frac{C(\theta, \epsilon)}{q^{2+\epsilon}},\]
      for all rational numbers $\frac{p}{q}$.
    \end{thm}

 \vspace{1cm}

 In higher variables, here is a famous theorem of Margulis, proving the Oppenhiem conjecture:

 \vspace{1cm}
 
 \begin{thm}
   If $Q$ is an indefinite and irrational quadratic form on $\R^n$, for $n \geq 3$, then $Q(\Z^n)$ is dense in $\R$.
 \end{thm}

 \vspace{1cm}

 There are no precise expectations for forms of  higher degrees!

 \newpage

\section{Subgroups of $\R^n$}
Recall that a finitely generated subgroup ${\mathcal L} \subset \R^n$ is a discrete subgroup if and only if the following equivalent conditions
are satisfied:

\begin{enumerate}

\item the mapping ${\mathcal L} \otimes \R \rightarrow \R^n$ is injective.
\item For any $\R$-subspace $V$ of $\R^n$ of dimension $d \leq n$,  ${\mathcal L}
  \cap V$ is generated by no more than
  $d$ elements.

\end{enumerate}

\vspace{1cm}

Further, we have the following assertions on subgroups $\LL$ of $\R^n$.

\begin{enumerate}

\item If $\LL$ is a closed subgroup of $\R^n$, then the connected component of identity,
  $\LL^0$ of $\LL$, is a subspace $\R^m \subset \R^n$, and $\LL/\LL^0$ is a discrete subgroup
  of $\R^n/\R^m$. Thus $\LL$ is uniquely described by giving its connected component
  $\LL^0$, and a discrete subgroup in $\R^n/\LL^0$. Therefore,  for a closed subgroup $\LL$ of $\R^n$,
  there is a decomposition of $\R^n$ as $\LL^0 \oplus \R^{n-m}$ with respect to which $\LL$
  is non-canonically of the form $\LL^0 \oplus {\mathcal M}$ for a discrete subgroup ${\mathcal M} \subset
  \R^{n-m}$.  

\item If $\LL$ is a general subgroup of $\R^n$ with $\bar{\LL}$ its closure, and  $\bar{\LL}^0 \cong \R^m$,
  its connected component of identity, then $\LL \cap \bar{\LL}^0$ is dense in $\bar{\LL}^0$,
  and $\LL= (\LL \cap \bar{\LL}^0) \oplus {\mathcal M}$ for a discrete subgroup ${\mathcal M} \subset
  \R^{n-m}$.  
\end{enumerate}

\newpage

    \section{Commutative algebraic groups}
    Connected commutative algebraic groups $\A$ over $\Q$ are built from the following three kinds of groups.

    \begin{enumerate}
    \item Finite dimensional vector spaces over $\Q$.
    \item Tori over $\Q$, which correspond bijectively (categorical equivalence)
      with finite free $\Z$-modules of finite rank equipped with an action of
      $\Gal(\bar{\Q}/\Q)$.
      \item Abelian varieties over $\Q$.
      
      \end{enumerate}

    The case of vector spaces over $\Q$ amounts to the study of lattices in $\R^n$ as in
    the last section,
    and there is nothing non-obvious to be said, so in the rest of the lecture, we will
    either deal with a torus or with an abelian variety. (I have not thought about a
    product of these, or the extension of one by the other.)

    In each case, we will consider a finitely generated (abelian!) subgroup $\mathcal L \subset \A(\Q)$,
    and consider the topological closure  $\bar{\mathcal L}$ of $\mathcal L$ in either
    $\A(\R)$ or $\A(\Q_p)$, and the answer that we desire is about $\bar{\mathcal L}^0$,
    which stands for

    \begin{enumerate}
    \item    the identity component of $\bar{\mathcal L}$  in $\A(\R) $,
    \item  the maximal compact subgroup of  $\bar{\mathcal L}$ in $\A(\Q_p)$.
    \end{enumerate}

    \vspace{1cm}

    The desired answer is that $\bar{\mathcal L}^0$ is:

    \begin{enumerate}
    \item a subgroup of    finite index in $\B(\R)$ where $\B \subset \A$ is an algebraic
      subgroup defined over $\Q$, and
    \item in the case of $\Q_p$, 
      we assume that  $\mathcal L \subset \A(\Q) $,
      hence $\mathcal L \subset \A(\Q_p)$, where it is contained in a compact subgroup.
      Then the desired answer is that up to commensurability,  $\bar{\mathcal L}^0 =  \Z_p^\ell$ where $\ell$ is
      the rank of ${\mathcal L}$, unless for some connected algebraic subgroup
      $\B \subset \A$,  ${\mathcal L} \cap \B(\Q_p)$ has rank which is greater than $\dim \B$.
          (Note that the maximal compact subgroup of $\B(\Q_p)$
    is up to commensurability, $\Z_p^d$, where $\dim \B = d$.)
    \end{enumerate}

\begin{conj} \label{reals}
  The desired answer holds if $\A$ is simple over $\Q$, or if $\A = R_{K/\Q}(\Bbb G_m)$.
  In fact if $\A$ is simple, it has no connected subgroups, and therefore for $\A$ simple,
  the conjecture is that:
  \begin{enumerate}
  \item In the case of $\A(\R)$, either ${\mathcal L} \subset \A(\Q)$ is a discrete subgroup
  of $\A(\R)$  (not an option if $\A$ is an abelian variety), or ${\mathcal L}$ is dense in $\A(\R)$.
\item In the case of $\A(\Q_p)$, assuming ${\mathcal L}$ contained in a maximal compact subgroup of
  $ \A(\Q_p)$, we have $\bar{\mathcal L} = \Z_p^\ell$ where $\ell$ is
      the minimum of the rank of ${\mathcal L}$ and $\dim \A$.
  \end{enumerate}
  \end{conj}

\begin{remark}
  The conjecture above in the case of  tori
  asserts that if certain tuple of algebraic numbers are mutiplicatively independent, then they are
  (equivalently their log's are) independent over $\R$, if there is enough room for $\R$-independence.  One might think
  that Schanuel's very general conjecture about algebraic independence of logarithms of algebraic numbers will
  solve any such question, such as Conjecture \ref{reals} above,  but that does not seem to be the case so far!
  \end{remark}

\begin{remark}
  In an earlier work, cf. \cite{Pr}, I looked at the corresponding question for general tori $T$ over $\Q$ for which $T(\R)/T(\Z)$ is compact
  and I looked at only those finitely generated subgroups $\Lambda$ containing $T(\Z)$, and proved analogous ``algebraicity theorem'' for 
  $\bar{\Lambda} \subset T(\R)$ assuming the Schanuel's conjecture. There is also the work of Waldschmidt in \cite{W2} characterizing
  those finitely generated subgroups $\Lambda$ of algebraic points of $T(\R)$ for which $\bar{\Lambda} = T(\R)$ proving an unconditional theorem. 
  \end{remark}

\begin{remark}
   Damien Roy, cf. \cite{Roy},  has made explicit construction of subgroups of $K^\times$ of rank $r_1+r_2+1$
  which are dense in $K_\infty^\times$. This paper also contains almost optimal bounds on the minimal
  rank of finitely
  generated subgroups of $K^\times$ which are dense in $K_S^\times$ for any finite set $S$ of places of $K$.
  \end{remark}

\section{Four exponential conjecture}

    Observe that the embedding  ${\Q}^\times \hookrightarrow \R^\times$,
  gives rise to an embedding:
  \[ ({\Q}^\times)^d   \hookrightarrow (\R^\times)^d,\]
  for all integers $d \geq 1$. If $d=2$, and ${\mathcal L}= \langle a_1, a_2\rangle$ with $a_1= (a_{11}, a_{12}) \in  ({\Q}^\times)^2,$
  and $a_2= (a_{21}, a_{22}) \in  ({\Q}^\times)^2,$ then for the matrix $M$,
 obtained using $ ({\Q}^\times)^2   \hookrightarrow (\R^\times)^2 \stackrel{\log}{\rightarrow} \R ^2,$
  $$M = \left ( \begin{array}{cc} 
  \log a_{11}  & \log a_{12}      \\
\log a_{21}  & \log a_{22}    
\end{array}
\right ),$$

\begin{enumerate}
\item the rows of $M$ are $\Q$ linearly independent if and only if the group ${\mathcal L}= \langle a_1, a_2\rangle$ is not generated
  by one element.

\item the columns of $M$ are $\Q$ linearly independent if and only if no nontrivial (algebraic) character of $\Q^\times \times \Q^\times$
  which is $\chi_{m,n}: (x,y)\rightarrow
  x^m\cdot y^n$ for $(m,n) \in \Z^2$ is trivial on  ${\mathcal L}= \langle a_1, a_2\rangle$.
\end{enumerate}

These are the conditions which are imposed in the famous ``four exponential conjecture'' (usually attributed to C. L. Siegel and T. Schneider)
and which then asserts that
the rank of the matrix $M$ is 2, equivalently, $\det M \not = 0$. Thus our conjecture (for ${\Bbb G}_m \times {\Bbb G}_m$ instead of  $ R_{K/\Q}({\Bbb G}_m)$)
is equivalent to the four exponential conjecture. The question being considered is simply if two vectors in $\R^2$, constructed using log of algebraic numbers
are $\R$-independent if they are $\Q$-independent?

\newpage

Just to give a particularly striking example, it is not known that the determinant of the following matrix (where $p,q$ are distinct primes $\geq 5$) is nonzero:

$$M = \left ( \begin{array}{cc} 
  \log 2  & \log 3      \\
\log p  & \log q    
\end{array}
\right );$$

or, in another form, if $2^t$ and $3^t$ both belong to $\Z$, then $t \in \Z$.

We refer to the book of Waldschmidt, cf. \cite{W1},  for an exposition on the four exponential conjecture, and the
six exponential theorem to which we come to now.

\section{The six exponential theorem}
A weaker version of the four exponential conjecture, is a theorem,
called the  ``six exponential theorem'' due to S. Lang and K. Ramachandra, and asserts that if
${\mathcal L}= \langle a_1, a_2, a_3\rangle$ with
$a_1= (a_{11}, a_{12}), a_2= (a_{21}, a_{22}), a_3= (a_{31}, a_{32}) \in  ({\Q}^\times)^2,$ then for the matrix $M$,
$$M = \left ( \begin{array}{cc} 
  \log a_{11}  & \log a_{12}      \\
  \log a_{21}  & \log a_{22}\\
  \log a_{31}  & \log a_{32}
  
\end{array}
\right ),$$
if the rows and the columns are $\Q$-independent, then the rank of $M$  is 2,   equivalently, the determinant of some $2 \times 2$
minor of $M$ is nonzero. Here too, the rows being $\Q$-linearly independent is equivalent to rank of the group  ${\mathcal L}= \langle a_1, a_2, a_3\rangle$ 
is equal to 3, and columns being   $\Q$-linearly independent is equivalent to saying no nontrivial (algebraic) character of $\Q^\times \times \Q^\times$
is trivial on   ${\mathcal L}= \langle a_1, a_2, a_3\rangle$, thus our conjecture (for ${\Bbb G}_m \times {\Bbb G}_m$ instead of  $ R_{K/\Q}({\Bbb G}_m)$)
for  ${\mathcal L}= \langle a_1, a_2, a_3\rangle$ is equivalent to the six exponential theorem.

\newpage

\section{Conjecture \ref{reals} for general tori}

Conjecture \ref{reals} is made for the torus $ R_{K/\Q}({\Bbb G}_m)$, or a simple torus,  and not a general
  torus over $\Q$ (though it could have been made) as the analogue of our conjecture for $({\Bbb G}_m)^d $ instead of  $ R_{K/\Q}({\Bbb G}_m)$
is not valid for $d \geq 3$, and the particular tori  that we have chosen is in the hope that for this
class of tori, there is no similar counter example!

To construct the counter-example, take $S= {\Bbb G}^3_m$, the 3 dimensional 
split torus over  ${\Bbb Q}$. The principle behind the counter-example is
the well-known observation that although the determinant of a 
 skew-symmetric $n \times n$ matrix 
consisting of  logarithm of algebraic numbers is 0 if $n$ is odd,
the rows and columns could be linearly independent over ${\Bbb Q}$, 
such as for the matrix:
$$
 \left(
\begin{array}{ccc}
0 & \log 2 & \log 3 \\
-\log 2 & 0 & \log 5 \\
-\log 3 & -\log 5 & 0 
\end{array}
\right).$$

For this matrix, the rows and columns are $\Q$-linearly independent, but the rank of this matrix is 2.
The subgroup ${\mathcal L}= \langle a_1, a_2, a_3\rangle$,
$a_1=(1,2,3), a_2=(1/2,1,5), a_3 =(1/3, 1/5, 1)$ is of rank 3 inside $\Q^{\times 3}$ and no nontrivial character of
 $S= {\Bbb G}^3_m$, is trivial on it.

\section{Some comments on the $p$-adic case: the Leopoldt conjecture}

\begin{remark}
  Conjecture 1 is the famous Leopoldt conjecture when ${\mathcal L} \subset K^\times$ is ${\mathcal L}= O_K^\times$.
  We recall that an equivalent form of the Leopoldt conjecture is the non-vanishing of the $p$-adic regulator which plays an
  important role in $p$-adic $L$-functions and Iwasawa theory. The Leopoldt conjecture is known basically only when $K/\Q$ is an abelian extension
  in which case the $p$-adic regulator is a product of linear forms in $p$-adic logarithm of algebraic numbers which are nonzero by a
  $p$-adic analogue of Baker's theory of linear forms in logarithm. However, even for abelian $K/\Q$, such a proof of Leopoldt conjecture does not generalise
  to prove the conjecture above to more general  subgroup ${\mathcal L} \subset K^\times$ unless the subgroup ${\mathcal L} \subset K^\times$
  has a subgroup of finite index which is freely generated by $\langle \ell^\sigma | \sigma \in \Gal(K/\Q)\rangle$ for some $\ell \in {\mathcal L}$.
    \end{remark}

\newpage
\section{Abelian varieties: a conjecture of Mazur}

We make some comments in the case of abelian varieties $A$ over $\Q$, for which by the Mordell-Weil theorem, $A(\Q)$ is a finitely generated
abelian group. Taking again a finitely generated subgroup ${\mathcal L} \subset A(\Q)$,
we consider the topological closure $\bar{\mathcal L}$ inside  $A(\Q_p)$. Recall that $A(\Q_p)$ is a compact abelian group
which is commensurable  
to $ \Z_p^n$ where $\dim A = n$, and therefore $\bar{\mathcal L}$ contains a subgroup of finite index which is
$\Z_p^d$ for some $d \leq n$. 
The question that interests us is when we can hope that the $\Z_p$-rank of (the free part of) $\bar{\mathcal L}$ is the same as the $\Z$-rank of $\mathcal L$.

\begin{conj} Let $A$ be a simple abelian variety over $\Q$ of $\dim n$,
  and ${\mathcal L} \subset A(\Q)$
  a finitely generated subgroup 
  of rank $m \leq n$.
  Then the topological closure $\bar{\mathcal L}$
  inside  $A(\Q_p)$ is isomorphic to
  $\Z_p^m$ (modulo finite torsion).
  \end{conj}

Here is the analogue of the previous conjecture,
now for $A(\Q)$ inside $A(\R)$ due to Mazur (\cite{M1}, \cite{M2}) who 
did not consider the $p$-adic case since  $\Z$ cannot be dense in $(\Z_p)^d$ for any $d > 1$, whereas $\Z$
can be dense in $(\ST^1)^d$ for any $d \geq 1$. Mazur only considers,
just as the group of full units in the Leopoldt conjecture,
the full Mordell-Weil group $A(\Q)$. Thus the following form of the conjecture of Mazur
is stronger than his.

\begin{conj} \label{mazur}(Strong form of Mazur's conjecture for Abelian varieties): Let $A$ be  a simple 
  abelian variety over ${\Q}$ of $\dim \geq 1$, then any non-torsion
element in $A(\Q)$ generates a subgroup of $A(\R)$ whose closure contains the connected component of identity
of $A(\R)$.
\end{conj}

In Conjecture 2 above, we could not have taken $A$ to be non-simple as the following example shows.
We will use the abelian variety $A=E\times E \times E$ of dimension 3, where $E$ is an elliptic curve containing
$\Z$-linearly independent
points $P,Q,R \in E(\Q)$, thus with the Mordell-Weil rank of $E$ being $\geq 3$. Let $S_1,S_2,S_3$ be rational points on $A(\Q) = E(\Q)\times E(\Q) \times E(\Q)$ defined by,
\begin{enumerate}
\item $S_1 = (0, P,Q),$
\item $S_2 = (-P, 0,R),$
  \item $S_1 = (-Q, -R,0).$
  \end{enumerate}
Then  the points $S_1,S_2,S_3$ are clearly $\Z$-linearly independent on $A$ but these are not
$\Z_p$-independent in $A(\Q_p)$. This one can prove by using a surjective group homomorphism with finite fibers
from $A(\Q_p)$ to $(\Z_p)^3$,  and noting that the images of $S_1,S_2,S_3$ must be $\Z_p$-dependent by the same
argument used earlier, that the determinant of a  $3\times 3$ skew symmetric matrix is zero! I am not sure how
to construct a counter-example over $\R$ of Conjecture \ref{mazur} for non-simple $A$.
\newpage

\section{A question of Madhav Nori}

Let $\pi$ be a finitely generated subgroup of $\GL_n(\Q)$. Then there is a
finite set $S$ of primes such that if $R=S^{-1}\Z$, then $\pi \subset \GL_n(R)$. 
There are two ways to measure of how large $\pi$ is:

\begin{enumerate}

\item the Zariski closure of $\pi$ in $\GL_n$ over $\Q$.
\item the topological closure $\widehat{\pi}$ of $\pi$ in $\GL_n(\widehat{R})$,
  where $\widehat{R}$ is the profinite completion of $R$.

\end{enumerate}

The theorems of Nori as well as those of Mathew-Vaserstein-Weisfeiler
imply that if $\G_{\pi}$ is the Zariski closure of $\pi$  in $\GL_n$ over $\Q$
which has finite fundamental group
then  the topological closure $\widehat{\pi}$ of $\pi$ in $\GL_n(\widehat{R})$
is open in $\G_\pi(\widehat{R})$.

On the other hand if $\G_\pi$ is a torus, then the above assertion is hopelessly false.
(For example, take a simple torus $T$ over $\Q$ of dimension 2. In this case,
the topological closure of the subgroup generated by a point of $T(\Q)$ cannot be
an open subgroup of $T(\Q_p)$.) However, Nori posed the following question, hoping it has an affirmative answer (and reducing it to the case of tori):

\vspace{2mm}

{\bf Question:} Is the dimension $d(p)$ of the $p$-adic Lie group obtained by taking the closure of $\pi$ in $\GL_n(\Z_p)$,  a constant function independent of $p$
where $p$ is a prime not belonging to $S$.

\vspace{2mm}

Nori noted the relationship of the above question to the Leopoldt conjecture.
It appears to me to be an excellent weaker version of Leopoldt conjecture,
valid in much greater generality, and valid also for Abelian varieties though
Nori was concerned only with linear algebraic groups.

However, Nori did not make precise
the number $d(p)$ to which we turn our attention to.

\section{The integer $d(p)$}

In this section, inspired by a remark of Madhav Nori, we make a conjecture
on $d(p)$ for which we turn to the notion of the {\it structural rank}
of a matrix due to Damien Roy, cf. \cite{Roy2}.

\begin{define}    Let $V$ and $W$ be finite dimensional vector spaces over $\Q$, with $V_\C=V \otimes \C$, and $W_\C=W \otimes \C$.
    Let $M \in \Hom_\C[V_\C,W_\C]$. Then the structural rank of $M$ is defined to be the smallest
    integer $s$ for which $M \in E \otimes \C$, where $E$ is a subspace of $\Hom_\Q[V,W]$
    consisting of homomorphisms of rank $\leq s$.
\end{define}

  \begin{remark}It is a theorem of Waldschmidt that if $r$ is the rank of the
    $d \times \ell$ matrix $M$ consisting of log of algebraic numbers, and $s$ is the structural rank, then
    \[ r \leq s \leq 2r,\]
    whereas conjecturally (such as by the Schanuel conjecture) $r=s$.
        \end{remark}

  Let $\mathcal{A}$ be a commutative algebraic group over $\Q$ with Lie
  algebra $\mathfrak{a}$ over $\Q$. Let,
  \[ {\rm exp} : \mathfrak{a} \otimes \R \longrightarrow \mathcal{A}(\R),\]
  be the exponential map.

  Let $\mathcal{L}$ be a finitely generated torsion free abelian group of rank $m$
  contained in  $\mathcal{A}(\Q)$. Then we have maps,
  
$$
  \xymatrix{
    & & \mathfrak{a} \otimes \R \ar[d]_-{\exp}  \\
\mathcal{L}=\Z^m \ar[r]^-{} & \mathcal{A}(\Q)   \ar[r]^-{} & \mathcal{A}(\R), }
$$ which allows us to define a map, call it $\log \mathcal{L}$ making the following diagram commute:

$$\xymatrix{
\mathcal{L}=\Z^m \ar[r]^-{\log \mathcal{L}} 
\ar[rd]^-{} & \mathfrak{a} \otimes \R \ar[d]_-{\exp}  \\
& \mathcal{A}(\R).
}
$$

  The following conjecture makes precise the common value of $d(p)$ in the question of Nori.

\begin{conj} \label {Nori}
  Let $\mathcal{A}$ be a commutative algebraic group over $\Q$, which is either a torus or is an abelian variety, with Lie
  algebra $\mathfrak{a}$ over $\Q$ with
  $ {\rm exp} : \mathfrak{a} \otimes \R \longrightarrow \mathcal{A}(\R),$
  the exponential map.
  Let $\mathcal{L}$ be a finitely generated torsion free abelian group of rank $m$
  contained in  $\mathcal{A}(\Q)$.
  For any prime $p< \infty$ such that the topological closure
  $\bar{\mathcal{L}}_p$ of $\mathcal{L}$
  inside  $\mathcal{A}(\Q_p)$ is compact, containing $\Z_p^{d(p)}$ as a subgroup of
  finite index, then such $d(p)$ are independent of $p$, and equal to the
 structural rank  of the linear map
    \[ \xymatrix{
      \Z^m \ar[r]^-{\log \mathcal{L}} &  \mathfrak{a} \otimes \R}. \]
\end{conj}

The above conjecture says nothing about the primes $p<\infty$
in $\Q$ where 
the topological closure
$\bar{\mathcal{L}}_p$ of $\mathcal{L}$ is non-compact. Nor does it say anything about  $\bar{\mathcal{L}}_p$ for $p = \infty$, to which we turn our attention. For this, we make the following definition.

  \begin{define} (Topological rank of a locally compact abelian group)
    Let $\mathcal{T}$ be a locally compact abelian group.
        Define the topological rank of   $\mathcal{T}$  to be the minimal
        rank of a free abelian subgroup  of  $\mathcal{T}$ whose closure in $\mathcal{T}$ is an open subgroup of finite index in  $\mathcal{T}$. If there are no abelian subgroups of finite rank in $\mathcal{T}$  whose closure has finite index, then we define the topological rank  of $\mathcal{T}$ to be infinity. (We will exclude such groups from consideration,
        such as $\Q_p$.)
          \end{define}

  \begin{example}
    The topological rank of $\Z_p^d$ is $d$, that of $(\ST^1)^d$ equals 1, that of $\Q_p^\times$ and $\R^\times$
    equal to 2.
    \end{example}
  With the notation as above, and denoting the common value of
  $d(p)$ in Conjecture \ref{Nori} to be $d(\mathcal{A})$, 
  here is the more precise assertion about
   $\bar{\mathcal{L}}_p$, $p \leq \infty$.

\begin{conj}
  Let $\T$ be a  torus  over $\Q$.
    Let $\mathcal{L}$ be a finitely generated torsion free abelian group of rank $m$
  contained in  $\T(\Q)$.
  For any prime $p< \infty$,
  let  the topological closure  $\bar{\mathcal{L}}_p$
  of $\mathcal{L}$
   inside  $\T(\Q_p)$ be of topological rank $\ell(p)$, and let
the topological closure  $\bar{\mathcal{L}}_\infty$ of $\mathcal{L}$
inside  $\T(\R)$ be of topological rank $\ell(\infty)$.
  Then
    \[ \ell(p) \geq d(\mathcal{A})  {\rm ~~for~~ all~}  p \leq \infty {\rm ~~at~~which~} \T {\rm~~ is~~split}. \]

 \end{conj}

\section{Going further: non-commutative groups}

In this section we consider a  general  algebraic group $\G$ over 
${\Bbb Q}$. Suppose $\G(\R)$ (or, $\G(\Q_p)$) operates on a space $X$ which is either a topological space or an algebraic variety
over $\R$ (or,  $\Q_p)$ transitively. 
Then the most general question that one could ask --- in the spirit of this paper ---
is about the closure of orbits of finitely generated subgroups $F$ of $\G(\R)$ (or $\G(\Q_p)$) consisting
of algebraic elements, and the question would be if these closures of orbits are ``locally'' orbits of closed, connected algebraic subgroups
of $\G(\R)$ (or, $\G(\Q_p)$). This is much too general a question, so one must restrict either $X$ and/or the finitely generated
subgroup $F$, to have any chance of a reasonable answer.

Recall that if $\G$ is a simple group over $\R$ or $\Q_p$, then by the Borel
density theorem, a lattice (i.e., a finite covolume discrete subgroup) in $\G(\R)$ or $\G(\Q_p)$ as the case may be, 
is Zariski dense in $\G$. As a consequence if $\LL$ is any finitely generated subgroup of $\G(\R)$ or $\G(\Q_p)$ as the case may be,
containing a lattice $\Lambda$, then the topological closure $\bar{\LL}$  in $\G(\R)$ or $\G(\Q_p)$ as the case may be,
is either discrete, i.e.,  $\bar{\LL} = \LL$ and contains $\Lambda$ as a subgroup of finite index, or is an open subgroup of
finite index in  $\G(\R)$ or $\G(\Q_p)$ as the case may be. (This follows by looking at the Lie algebra of $\bar{\LL}$ which 
is a module for the Zariski dense subgroup  $\Lambda$, hence for the Lie algebra of $\G$.  Therefore, $\G$ being simple,
the Lie  algebra of $\bar{\LL}$ must be either trivial
or that of $\G$.)

\vspace{1cm}

Thus considering topological closure of finitely generated subgroups of $\G(\R)$ or $\G(\Q_p)$ containing
lattices is not so interesting. However, instead of considering finitely generated subgroups of  $\G(\R)$ or $\G(\Q_p)$ 
containing a lattice $\Lambda$, we could consider   $\G(\R)/\Lambda$ or $\G(\Q_p)/\Lambda$, and look
at the topological closure of the image of finitely generated subgroups $\LL$ of $\G(\R)$ (resp., $\G(\Q_p)$)
in  $\G(\R)/\Lambda$ (resp., $\G(\Q_p)/\Lambda$); this is the same as closure of orbits of $\LL$ on
$\G(\R)/\Lambda$ or $\G(\Q_p)/\Lambda$.  In line with 
the theme of this paper, we should be considering
the topological closure of the image of finitely generated subgroups consisting of algebraic points
of $\G(\R)$ (resp. $\G(\Q_p)$), and the hope would be that the closure of the image of $\LL$ in  $\G(\R)/\Lambda$ or $\G(\Q_p)/\Lambda$
as the case may be, equals, up to finite index,
image of $\H(\R)$ or $\H(\Q_p)$ for $H$ an algebraic subgroup $\H \subset \G$ under the projection to  $\G(\R)/\Lambda$ or $\G(\Q_p)/\Lambda$.
Note that the  topological closure of $\LL$ in  $\G(\R)/\Lambda$, is the image in  $\G(\R)/\Lambda$
of the
closure of the set $\LL \cdot \Lambda$ in $\G(\R)$. Similarly,
the topological closure of the orbit of $\LL$ passing through the point $x \Lambda$ in  $\G(\R)/\Lambda$ is the image  in  $\G(\R)/\Lambda$
of the topological closure of  $\LL \cdot x \cdot  \Lambda$ in $\G(\R)$. Thus the question from this perspective is to understand
closure of  $\LL \cdot x \cdot  \Lambda = x \cdot (x^{-1}\LL \cdot x) \cdot  \Lambda $ in $\G(\R)$ ($x$ an algebraic element of $\G(\R)$)
where $\LL$ is any finitely generated subgroup of $\G(\R)$ consisting of algebraic
elements which after replacing $\LL$ by  $x^{-1}\LL \cdot x$ is equivalent to understand the closure of $\LL \cdot \Lambda$ where $\LL$ is a finitely
generated subgroup of $\G(\R)$ consisting of algebraic elements, and $\Lambda$ is a lattice in $\G(\R)$.

\vspace{1cm}

We would like to suggest the 
analogue of our conjectures for commutative algebraic groups to  assert that
for $\G$ a semisimple algebraic group over $\Q$, $\Lambda$ an arithmetic subgroup
of $\G(\Q)$, the
closure of the image 
in $ \G({\Bbb R})/\Lambda$ (resp., $ \G(\Q_p)/\Lambda$)
of a finitely generated subgroup 
$\LL$ of $\G(\R)$ (resp., $\G(\Q_p)$) consisting of algebraic elements
is of the form $\Gamma_\H \backslash \H({\Bbb R})$ (resp., $\Gamma_\H \backslash \H(\Q_p)$) 
for an algebraic subgroup $\H$ of
$\G$ defined over ${\Bbb Q}$. (We recall  a well-known theorem in Lie theory that $\H(\R) \cdot \Lambda$
is a closed subset of $\G(\R)$ if  $\H(\R) \cap \Lambda$ is a lattice in $\H(\R)$, or if $\H$ is defined over
$\Q$ and $\Lambda$ is an arithmetic group in $\G(\Q)$.)

We remark that our suggested analogue contains a consequence of
M. Ratner's theorem (the proof of the so-called Raghunathan conjecture)
as observed by Dani and Raghunathan, cf. Cor. 4.9 in \cite{V} in a very special
case. It states that if a semi-simple group G over ${\Bbb R}$ (with 
$\G({\Bbb R})$ non-compact) has two
distinct ${\Bbb Q}$ structures, with corresponding lattices
${\Gamma_1}$ and ${\Gamma_2}$, then $\Gamma_1 \cdot \Gamma_2$ is 
dense in $\G({\Bbb R})$ (in the Euclidean topology); similarly for $\G(\Q_p)$ in place of $\G(\R)$.

\end{document}